\documentclass[12pt]{article}
\usepackage[letterpaper,margin=1in]{geometry}
\usepackage[T1]{fontenc}
\IfFileExists{newtxtext.sty}{\usepackage{newtxtext,newtxmath}}{\usepackage{mathptmx}}
\usepackage{setspace}
\usepackage{microtype}
\newcommand{\fasection}[1]{%
  \vspace{1.4em}
  {\centering\textbf{#1}\par}
  \vspace{0.6em}\noindent}

\begin{document}

\begin{center}
  {\LARGE\bfseries Automation Without Understanding}\\[0.8em]
  {\large\itshape Why the United States must preserve mathematical
   capacity in the age of AI}\\[1.2em]
  {\normalsize\scshape Jun--Yong Park}



\end{center}

\vspace{1.2em}

\noindent No one thinks about clean water until it turns brown. No one thinks about functioning courts until a verdict is corrupt. The defining feature of infrastructure is its invisibility: when it works, nobody notices; when it fails, everything downstream fails with it. Mathematical understanding is infrastructure of this kind. It is invisible and integral. And it is being dismantled, not by malice but by neglect, at the precise moment when the AI systems built on top of it are becoming more consequential and more opaque than any technology in history.

\smallskip

Two developments are unfolding at once. Artificial intelligence systems have begun to do genuine mathematics: not calculation but discovery, of a kind that until recently only trained human mathematicians could produce. And the United States, through budget chaos and institutional drift, is degrading the pipeline that produces humans capable of understanding what the machines are doing. Either development alone would deserve attention. The second, unfolding alongside the first, amounts to a strategic error.

\fasection{THE MACHINES CAN DO MATHEMATICS}

On May 20 of this year, OpenAI announced a breakthrough on the planar unit distance problem in discrete geometry~\cite{openai2026}, first posed by the Hungarian mathematician Paul Erd\H{o}s in 1946~\cite{erdos1946}. The question sounds simple: scatter $n$ points on a plane and ask how many pairs of them can sit exactly one unit apart. Erd\H{o}s conjectured a ceiling on that count, and for eighty years the best minds in combinatorics believed him. The model proved him wrong. It constructed families of point arrangements that break the conjectured ceiling, and it did so by importing machinery from algebraic number theory, including class field towers, a corner of mathematics built for entirely different purposes.

\smallskip

The claim did not rest on the company's word. The same day, nine mathematicians published a companion paper~\cite{companion2026,sawin2026} verifying the argument and translating it into human mathematical language. Among them was the Fields Medalist Timothy Gowers, who wrote that had a human submitted the result to the \textit{Annals of Mathematics}, the field's most prestigious journal, he would have recommended acceptance without hesitation.

\smallskip

Nor was the result a stunt. The system in question was not designed for mathematics; it was a general-purpose reasoning model powerful enough to do what no human had done. And the trajectory has been steep. In 2024, Google DeepMind's specialized systems reached the silver-medal standard on International Mathematical Olympiad problems~\cite{deepmind2024}. A year later, general-purpose models from both OpenAI and Google DeepMind cleared the gold-medal standard~\cite{deepmind2025,openaigold2025}. Today, frontier models produce non-obvious constructions, identify connections across fields, and contribute to research-level work.

\smallskip

The machines are doing mathematics. That sentence was controversial three years ago. It is simply true today. The question is no longer whether they can. The question is whether the United States will retain the human capacity to understand, verify, direct, and, when necessary, refuse what the machines produce.

\fasection{THE QUIET DISMANTLING}

At exactly this moment, that capacity is being run down. Last year, the administration proposed cutting the National Science Foundation's budget by more than half; within that request, funding for the mathematical and physical sciences would have fallen by roughly two-thirds, and for the mathematical sciences specifically, analysts put the reduction above 70 percent. Congress ultimately rejected most of the cuts: the appropriation signed in January keeps the NSF roughly flat, at \$8.75 billion against roughly \$9 billion the year before.

\smallskip

But the damage did not wait for the final number. During the months of uncertainty, more than \$14 million in grants already promised to mathematics programs was clawed back, according to an analysis by \textit{Scientific American}~\cite{hasson2025}. The NSF canceled funding for a national research symposium for women in mathematics four business days before it was scheduled to begin. The American Mathematical Society, a professional society rather than a funding agency, had to stand up \$1 million in emergency ``backstop'' grants to keep stranded programs alive. And universities, facing the same pressure from every federal direction at once, pulled back: Harvard, Princeton, MIT, and the entire University of California system imposed hiring freezes.

\smallskip

Flat funding after a year of institutional chaos is not recovery. It is stabilization at a lower level of confidence, and even the stabilization has proven illusory. By June, according to reporting in \textit{Science}, the NSF was quietly cutting hundreds of its basic research programs, in what agency insiders suspected was a bid to free more than a billion dollars for a new commercialization initiative that Congress had never funded~\cite{mervis2026,nsfxlabs2026}. Some fields were cut far beyond the 5 percent limit that lawmakers had written into the appropriation itself, and some programs did not learn their actual budgets until five months after the money was approved. And the destination of the money completes the irony: the new initiative exists to commercialize quantum systems and AI-driven scientific instrumentation. The government is funding the machines with money taken from the people who understand them.

\smallskip

The universities, meanwhile, have already restructured around the expectation of decline. Retirements go unreplaced. Doctoral cohorts shrink: George Washington University admitted zero new funded mathematics doctoral students this year~\cite{ihe2026,hatchet2026}, its funding sufficient only for the students it already had, while Harvard cut its doctoral intake by more than half and Chicago announced a 30 percent reduction~\cite{chronicle2025}. The pipeline that produces the next generation of mathematical minds is narrowing just as the machines arrive. 

\smallskip

None of this reflects a deliberate policy aimed at mathematics. It reflects a convergence of budget cuts, institutional restructuring, and a reasonable-sounding inference: if the machines can do mathematics, why pay humans to do it?

\smallskip

The answer is the same one that keeps engineers who understand water chemistry inside an automated treatment plant, and human jurists on the bench even if machines could draft flawless legal opinions. Automation works until it does not. Every sophisticated system eventually fails in a way its designers did not foresee, and when it does, the country needs people who understand the system from the inside: people who can trace the failure down to first principles. In domains of judgment, the need is not only technical but civic. A verdict that no human understands is a verdict no institution can stand behind. A system that nobody understands is a system that nobody can govern.

\fasection{A PROOF IS A PRODUCT. UNDERSTANDING IS A CAPACITY.}

Mathematical capacity is not a stockpile of theorems. It is a community of people transformed, by years of rigorous training, into minds able to see abstract structure inside complexity, detect errors that formal systems miss, and reason with unusual precision. No one learns to think as a mathematician by reading a book. One learns it by spending years confronting theoretical frameworks initially beyond one's reach, under the supervision of people who learned the same way, inside institutions that believe such training is worth sustaining.

\smallskip

Remove the funding and weaken the institutions, and the training pipeline breaks. When it breaks, the capacity is not preserved elsewhere; it becomes thinner, more uneven, and far harder to rebuild. A country cannot conjure a mathematical workforce on demand any more than it can conjure an officer corps or a cadre of nuclear engineers on demand. The training takes years. The intellectual tradition takes generations.

\smallskip

The usual reply is that future students can simply learn from recorded lectures, archived textbooks, or AI tutors. But mathematical capacity is not a body of information; it is a practiced form of judgment. Every surgical procedure ever performed could be captured on video. No one would attempt open-heart surgery on that basis.

\smallskip

A proof is a product. Understanding is a capacity. Machines may increasingly manufacture the product. The capacity to interpret, test, extend, and challenge it must still be cultivated in human minds and sustained by institutions that train those minds. When the people who carry that judgment are gone, the books remain on the shelf. The capacity to read them does not.

\smallskip

The Erd\H{o}s episode itself makes the point. The companion paper's authors were explicit that human experts had to discuss, simplify, and improve the machine's argument before the field could absorb it. The proof became knowledge, something other mathematicians can now build on, only because a community existed with the training to receive it. That community is the asset, and every additional machine-generated theorem raises its value.

\smallskip

The training produces habits that matter far beyond mathematics. A mathematician learns to ask questions that are often unwelcome in fast-moving institutions: What, exactly, is the object under discussion? Which assumptions are doing the work? What would count as a counterexample? Does the conclusion actually follow from the premises? Is the apparent pattern structural, or merely accidental? These are precisely the habits required to audit a model, expose a hidden assumption, distinguish a robust result from an artifact of the data, and recognize when a formal system is being asked a question it was not built to answer.

\smallskip

Pure mathematics also preserves a culture of unusually serious criticism. In a good seminar, authority is not enough. A famous speaker and a first-year graduate student are answerable to the same demands: define the terms, justify the inference, make the argument clear. That culture cannot be recreated by buying more chips or allocating more capital. It requires time, institutional continuity, and a community large enough to sustain standards that no individual can maintain alone.

\smallskip

That is what is at risk.

\fasection{EMERGENCE WITHOUT UNDERSTANDING}

The training objective behind today's large language models is deceptively simple: predict the next token in a sequence. Yet systems trained at sufficient scale, then refined through post-training and extended inference, display forms of abstraction, planning, and multistep problem-solving that their creators do not fully understand. They never stop predicting tokens, but the behavior that emerges is increasingly difficult to describe in the vocabulary that built it. A system trained to continue text can now overturn an eighty-year-old conjecture. The gap between the simplicity of the rule and the complexity of the result is one of the central scientific facts of the present moment.

\smallskip

Serious work is under way on interpretability, formal verification, and the theory of learned representations. But these efforts remain immature relative to the systems they are trying to explain. There is not yet a systematic mathematical theory of how large artificial neural networks form internal concepts, why certain capabilities appear abruptly, or how to predict which forms of reasoning will emerge at a given scale.

\smallskip

There are glimpses of what such a theory would have to explain. A few years ago, researchers trained a small artificial neural network on modular arithmetic, a staple of undergraduate number theory, and then reverse-engineered what it had learned~\cite{nanda2023}. The network had not memorized its answers. It had independently converged on an algorithm built from discrete Fourier transforms and trigonometric identities, converting addition into rotation about a circle: mathematical structure that no one had taught it and that its designers did not anticipate. The episode is celebrated in the young field of mechanistic interpretability, but its most important lesson is usually missed. Decoding the network was possible only because the researchers could recognize Fourier structure when they saw it. The machine's inner workings became legible only to humans who had the mathematics to read them. And that was a toy model; frontier systems are larger by many orders of magnitude. We cannot map what a frontier model has learned with the mathematical equivalent of a magnifying glass.

\smallskip

When mathematicians encounter an unfamiliar object, they do not ask only whether it performs well. They ask what kind of object it is: What are its invariants? Its symmetries? Its limiting behaviors? Which examples are canonical, and which apparent phenomena are artifacts of a particular presentation? Nothing like that framework yet exists for the internal representation spaces of frontier AI systems. Existing tools reveal fragments: circuits, features, scaling curves. But fragments are not a theory. We can inspect portions of the machine without knowing what sort of mathematical object the whole machine has become.

\smallskip

The institutional picture compounds the scientific one. Frontier models were trained on a vast public inheritance: papers, textbooks, lecture notes, and solved problems produced largely in universities and public research institutions. Yet the work of building and studying the most capable systems now occurs inside a small number of private firms, with limited outside access to the models, data, and experiments that independent scrutiny requires. The public intellectual infrastructure capable of understanding these systems has not kept pace with the private infrastructure capable of building them. Societies are coming to rely on systems whose capabilities are advancing faster than the science that explains them.

\smallskip

Pure mathematics cannot close this gap alone; neither can machine learning, computer science, statistics, formal methods, security engineering, or governance. But mathematics belongs at the center of the coalition, because the task is not only to make the systems more powerful. It is to develop the conceptual frameworks by which their structure can be understood, their claims audited, and their limits known.

\fasection{THE FIRST-CONTACT PROBLEM}

Artificial superintelligence may arrive soon, late, or never. No one knows, and a responsible argument does not pretend otherwise. The more immediate problem is already here: institutions are beginning to delegate judgment to systems whose reasoning they cannot fully examine. A model may advise on military logistics, financial risk, scientific research, cybersecurity, or critical infrastructure. It may be right often enough to become indispensable. Indispensability without intelligibility is a dangerous bargain.

\smallskip

When such a system produces a recommendation with serious consequences, ``the model said so'' cannot be the end of the discussion. Someone must be able to ask what was optimized, what evidence was omitted, which failure modes remain plausible, and whether the system answered the question it was actually asked. That is already difficult. It will become harder as the systems grow more capable.

\smallskip

Suppose, then, that artificial superintelligence does arrive. Americans should not imagine that human beings will meet such a system as intellectual equals, or that guardrails written by weaker minds will permanently constrain a mind that understands them better than their authors do. Nor is mathematics a magic diplomatic language through which humans can negotiate obedience. But neither does human understanding become irrelevant. Mathematics remains humanity's most precise language for structure, constraint, proof, and consequence: the language in which one can state exactly what follows from what, which assumptions are required, where an argument breaks, and what remains invariant under change. A society that loses the people who think in those terms loses part of its power to inspect minds it cannot otherwise comprehend.

\smallskip

The first-contact problem, in its sober form, is not a cinematic encounter with a machine god. It is the gradual normalization of consequential decisions whose real logic is inaccessible to the human beings expected to govern them. A model can produce an explanation in fluent English. But fluency is not transparency: natural language can make an argument sound persuasive without making it valid. The remedy, to which I will return, is to demand something better than persuasion.

\smallskip

No community of mathematicians will make humanity equal to an artificial superintelligence. But a society without such a community will approach increasingly powerful systems with less than it has now: fewer people able to recognize when a result is profound, when a proof is incomplete, when a formalization has smuggled in the wrong question, and when apparent intelligence masks a deeper failure of understanding. The danger is not simply that machines may become too capable. It is that human institutions may become too intellectually dependent to know what those capabilities are being used for.

\fasection{WHAT WASHINGTON SHOULD DO}

The United States should treat mathematical capacity as a strategic asset, on a par with semiconductor capability, national-security research, and energy security. That means more than preserving a floor of funding. Critical infrastructure is not run at minimum margins; it is built with redundancy, because future failures cannot be predicted.

\smallskip

First, preserve the pipeline that produces mathematical researchers: strong undergraduate programs, funded graduate training, postdoctoral positions, permanent academic jobs, and research institutes capable of supporting long-horizon work insulated from commercial pressure. A civilization that depends on mathematics but ceases to cultivate mathematicians is spending down an endowment it has stopped replenishing.

\smallskip

Second, build an independent national capacity for AI assurance: public institutions where mathematicians, statisticians, computer scientists, and verification researchers can evaluate frontier systems with their own access to models, data, and experiments. The people responsible for judging those systems cannot be wholly dependent on the companies building them.

\smallskip

Third, stop demanding that every mathematical project justify itself through immediate application. Some of the mathematics that becomes indispensable twenty years from now will look obscure today; that is not evidence of irrelevance, it is how foundational knowledge has always worked. The point is not to convert every pure mathematician into an AI engineer. It is to preserve the ecosystem from which new concepts, unexpected methods, rigorous critiques, and independent standards of proof emerge. The Erd\H{o}s proof itself ran through class field towers: machinery a century in the making, built by people who could not have named its eventual application.

\smallskip

Fourth, and most concretely: when AI systems are asked to perform consequential reasoning, they should be required to expose their decision-critical claims in a formal, machine-checkable form, not merely as persuasive natural-language explanations. This is the demand for something better than persuasion. A formally specified claim, accompanied by a proof or a verifiable certificate, can be checked by independent tools and human auditors without relying on the system's own account of why it is right. This would not solve the alignment problem: a valid proof can still rest on false premises, optimize the wrong objective, or faithfully answer a badly specified question. But it would change the nature of oversight, converting part of AI reasoning from opaque persuasion into auditable structure, and it would create an enduring role for independent experts trained in mathematics, formal verification, computer science, and the relevant domains.

\smallskip

Nor is that expertise incidental to the checking machinery itself. Formal verification does not replace human mathematical capacity; it consumes it. When one research team set an automated system the task of formalizing the Erd\H{o}s disproof in Lean, a leading proof-assistant language, the attempt collapsed in an instructive way~\cite{leanmarathon2026}: the necessary algebraic number theory had never been encoded in Lean's human-built libraries, and the automation papered over the gap with placeholder structures that type-checked while proving nothing. Another team~\cite{aleph2026} did complete a formalization, but only by taking the two deepest class-field-theoretic results as explicit unproven hypotheses, a hundred years of human work written directly into the theorem's signature as trust. Every formal library is centuries of human mathematics, translated line by line by people trained to know what the lines mean. Defund the training, and the verification infrastructure hollows out with it.

\smallskip

None of this would be unprecedented. In 1984, a National Research Council committee chaired by Edward David, a former presidential science adviser, warned in \textit{Renewing U.S. Mathematics: Critical Resource for the Future}~\cite{david1984} that federal support for the field had quietly eroded even as the country's dependence on it grew. Congress and the NSF responded, and funding rose substantially in the years that followed, though it never reached the report's targets and the commitment faded within a decade: a reminder that capacity is not restored by a single appropriation. The David report treated mathematical talent as a national resource, and Washington, for a time, acted accordingly.

\smallskip

Others have absorbed that lesson even as the United States forgets it. In 2019, four of China's top science agencies jointly issued a national plan to strengthen mathematical research~\cite{china2019}, including sustained support for pure mathematics, on the stated grounds that the discipline underpins aerospace, biomedicine, energy, and artificial intelligence. The five-year plan covering 2021 to 2025~\cite{fyp2021} went on to name mathematics first among the basic disciplines for which Beijing would build frontier research centers and train elite students. China does not fund mathematicians out of sentiment.

\smallskip

My father used to tell me: ``When you cannot tell what is the right thing to do, do the most difficult thing well. That is the right thing to do.'' Mathematical understanding is difficult. It is slow. It is unprofitable on any timescale a budget can see. It takes a kind of passion that no benchmark can measure. And it is part of the foundation on which everything else stands. 

\smallskip

The choice facing the United States is not between human thought and machine intelligence. It is between automation with understanding and automation without it. A country that chooses the second path does not become more advanced. It becomes more dependent, more fragile, and less capable of directing its own future.



 \vspace{+10pt}

\noindent Jun--Yong Park \enspace ---\enspace \texttt{june.park@sydney.edu.au} \\
\textsc{School of Mathematics and Statistics, University of Sydney, Australia}


\medskip

\begin{thebibliography}{99}

\bibitem{erdos1946}
P.~Erd\H{o}s, \textit{On Sets of Distances of $n$ Points}, 
The American Mathematical Monthly, \textbf{53}, No. 5, (1946): 248--250.

\bibitem{openai2026}
OpenAI, \textit{An OpenAI model has disproved a central conjecture in discrete geometry}, May 20, 2026. \texttt{openai.com/index/model-disproves-discrete-geometry-conjecture}.

\bibitem{companion2026}
N.~Alon, T.~F.~Bloom, W.~T.~Gowers, D.~Litt, W.~Sawin, A.~Shankar, J.~Tsimerman, V.~Wang, and M.~Matchett Wood, \textit{Remarks on the disproof of the unit distance conjecture}, arXiv:2605.20695, May 2026.

\bibitem{sawin2026}
W.~Sawin, \textit{An explicit lower bound for the unit distance problem}, arXiv:2605.20579, May 2026.

\bibitem{deepmind2024}
Google DeepMind, AlphaProof and AlphaGeometry teams, \textit{AI achieves silver-medal standard solving International Mathematical Olympiad problems}, July 2024. \\
\texttt{deepmind.google/blog/ai-solves-imo-problems-at-silver-medal-level}.

\bibitem{deepmind2025}
Google DeepMind, Thang Luong and Edward Lockhart, \textit{Advanced version of Gemini with Deep Think officially achieves gold-medal standard at the International Mathematical Olympiad}, July 2025. 

\bibitem{openaigold2025}
OpenAI, announcement of gold-medal performance on the 2025 International Mathematical Olympiad, July 2025.

\bibitem{hasson2025}
E.~R.~Hasson, \textit{Can U.S. math research survive NSF funding cuts?}, Scientific American, July 2025.

\bibitem{mervis2026}
J.~Mervis, \textit{Exclusive: NSF slashes research programs to support new tech initiative, insiders say}, Science, June 2026.

\bibitem{nsfxlabs2026}
National Science Foundation, announcement of the NSF X-Labs program, May 14, 2026.

\bibitem{ihe2026}
J.~Alonso, \textit{George Washington U pauses admissions to 5 Ph.D. programs}, Inside Higher Ed, January 26, 2026.

\bibitem{hatchet2026}
G.~Jakubowski, \textit{CCAS to shrink, halt doctoral program admissions in ``devastating'' 7 percent package cut}, The GW Hatchet, January 26, 2026.

\bibitem{chronicle2025}
M.~Zahneis, \textit{Has the graduate-school collapse begun?}, The Chronicle of Higher Education, December 3, 2025.

\bibitem{nanda2023}
N.~Nanda, L.~Chan, T.~Lieberum, J.~Smith, and J.~Steinhardt, \textit{Progress measures for grokking via mechanistic interpretability}, arXiv:2301.05217, 2023.

\bibitem{leanmarathon2026}
Y.~Zhang, Y.~Sun, T.~Suzuki, J.~D.~Lee, and F.~Liu, \textit{LeanMarathon: Toward reliable AI co-mathematicians through long-horizon Lean autoformalization}, arXiv:2606.05400, June 2026.

\bibitem{aleph2026}
A.~Fetisov, \textit{Aleph Prover Formalized Planar Unit Problem Disprove}, Logical Intelligence, May 28, 2026. Repository: \texttt{github.com/logical-intelligence/erdos-unit-distance}.

\bibitem{david1984}
National Research Council, Ad Hoc Committee on Resources for the Mathematical Sciences (E.~E.~David, Jr., chair), \textit{Renewing U.S.\ Mathematics: Critical Resource for the Future}, National Academy Press, Washington, D.C., 1984.

\bibitem{china2019}
Ministry of Science and Technology, Ministry of Education, Chinese Academy of Sciences, and National Natural Science Foundation of China, \textit{Work Plan for Strengthening Mathematical Science Research}, Guo Ke Ban Ji [2019] No.~61, July 2019. 

\bibitem{fyp2021}
\textit{Outline of the 14th Five-Year Plan (2021--2025) for National Economic and Social Development and Long-Range Objectives Through the Year 2035}, adopted by the National People's Congress, March 2021.

\end{thebibliography}
\end{document}